\documentclass[11pt]{article}
\usepackage{amssymb,latexsym,color,amsmath,pifont,epsfig,cite}
\usepackage[top=0.8in,bottom=1in,left=1in,right=1in]{geometry}
 \usepackage{graphicx}   
\usepackage{epsfig}
\usepackage{graphicx,graphics}
\usepackage{mathrsfs}
\usepackage{times}
\usepackage{mathptmx}
\usepackage{amssymb,amsmath,amsthm}
\usepackage{amsfonts}
\usepackage{txfonts}
\usepackage{subfigure}
\usepackage{amsfonts} \usepackage{txfonts}
\newtheorem{theorem}{\textbf {Theorem}}
\newtheorem{lemma}{\textbf {Lemma}}

\newtheorem{remark}{\textbf{Remark}}
\newtheorem{example}{\textbf{Example}}

\begin{document}
\title{On Spatially Uniform Behavior in Reaction-Diffusion PDE and Coupled ODE Systems}
\author{Murat Arcak\thanks{Department of Electrical Engineering and Computer Sciences,
University of California, Berkeley. Email:  arcak@eecs.berkeley.edu. Research supported in part by the National Science Foundation under grant ECCS 0852750 and by the Air Force Office of Scientific Research under grant FA9550-09-1-0092.}}
\date{\today}
\maketitle

\begin{abstract}
We present a condition which guarantees spatial uniformity for the asymptotic behavior of the solutions of a reaction-diffusion PDE with Neumann boundary conditions.  This condition makes use of the Jacobian matrix of the reaction terms and the second Neumann eigenvalue of the Laplacian operator on the given spatial domain, and replaces the global Lipschitz assumptions commonly used in the literature with a less restrictive Lyapunov inequality.  We then present numerical procedures for the verification of this Lyapunov inequality and illustrate them on models of several biochemical reaction networks.  Finally, we derive an analog of this PDE result for the synchronization of a network of identical ODE models  coupled by diffusion terms.
\end{abstract}

\section{Introduction}
Spatially distributed system models are essential for many fields of science and engineering.
In cell biology, gradients of protein activities organize signaling around cellular structures and provide positional cues for important processes, such as cell division \cite{Kholodenkodiffusion06}.  One of the theories for spatial organization and pattern formation is based on diffusion-driven instability \cite{turing, SegJac72}, which has been a subject of intense study as surveyed in \cite{murray89,CroHoh93,OthPaiUmuXue09}.  This phenomenon occurs when one of the higher spatial modes in the reaction-diffusion partial differential equation (PDE)
 is destabilized by diffusion, thus causing nonuniformities to  grow.  Understanding when the solutions of a reaction-diffusion PDE exhibit uniform behavior is an important problem, because it rules out diffusion-driven instabilities and justifies a simpler ordinary differential equation (ODE) modeling. The standard approach to proving spatial uniformity in the literature is to establish
 exponential decay of initial nonuniformities by using global Lipschitz bounds on the vector field representing reaction terms \cite{othmer77,AshOth78,ConHofSmo78,JonSle83}.

In the first part of this paper, we study the reaction-diffusion PDE:
\begin{equation}\label{rdnet}
\frac{\partial x}{\partial t}=f(x)+D\nabla^2x,
\end{equation} subject to
Neumann boundary conditions and other technical assumptions detailed in Section \ref{pdesec}, and give a condition for uniform behavior of the solutions that does not rely on a global Lipschitz assumption on $f(x)$. Instead, our main result (Theorem \ref{network}) requires that a Lyapunov inequality be satisfied by the matrix $J(x)-\lambda_2D$, where
\begin{equation}
J(x):=\frac{\partial f(x)}{\partial x}
\end{equation}
is the Jacobian and $\lambda_2$ is the second Neumann eigenvalue of the operator $L=-\nabla^2x$ on the given spatial domain.  Even when the global Lipschitz condition of \cite{othmer77,AshOth78,ConHofSmo78,JonSle83} holds, our result can achieve orders of magnitude improvements over the estimates obtained from this Lipschitz bound (see Example \ref{Goodwinex} for a comparison).

In the second part of the paper (Section \ref{LMIsec}), we parameterize $J(x)$ with constant matrices and develop procedures to verify the Lyapunov inequality employed in Theorem \ref{network}.
The first procedure, described in Theorem \ref{LMIlem}, incorporates $J(x)$ within convex and conic hulls of constant matrices and derives a {\it linear matrix inequality} (LMI) \cite{boyd} for the vertices.  The second procedure, presented in Theorem \ref{dyad}, studies a special convex set and reduces the dimension of the LMI in Theorem \ref{LMIlem}.  For reaction networks that exhibit special structures, the LMI in Theorem \ref{dyad} is also amenable to analytical feasibility tests.
One such test is illustrated in Example \ref{Goodwinex} on a variant \cite{Thr91} of Goodwin's model \cite{Goodwin65} for oscillations in enzyme synthesis.  In Example \ref{Goldbeterex}, we study a model by Goldbeter \cite{gold95} for circadian rhythms and investigate the feasibility of the LMI numerically. 

In a recent study \cite{JovArcSon08}, we gave conditions for the stability of the spatially uniform fixed point for reaction-diffusion systems where the reaction terms exhibit a cyclic structure.  In the present paper we do not restrict ourselves to cyclic reactions and, more importantly, we do not require that the attractor be a fixed point.  Indeed, the reactions in Examples \ref{Goodwinex} and \ref{Goldbeterex} exhibit limit cycles and Theorem \ref{network} guarantees spatial uniformity of the oscillations rather than stability of a fixed point.

In the third part of the paper (Section \ref{odesec}), we derive an analog of Theorem \ref{network} for a finite number of identical ODEs coupled via diffusion-like terms \cite{Hal97}.  For ODEs, the equivalent of spatially uniform behavior is {\it synchronization}, on which a large literature exists as reviewed in \cite{sync}.  Our main result (Theorem \ref{odethm}) in this part employs the same condition as Theorem \ref{network}, where $\lambda_2$ now represents the second smallest eigenvalue of the Laplacian matrix for the graph describing the coupling of the subsystems.  The proof of this result exploits properties of the Laplacian matrix that are analogous to those of the Laplacian operator employed in Theorem \ref{network}. 
In Example \ref{rodolphe}, we make a connection between Theorem \ref{odethm} and the {\it incremental passivity} approach to synchronization employed in \cite{StaSep07}.

\section{Uniform Behavior in Reaction-Diffusion PDEs}\label{pdesec}
Let $\Omega$ be a bounded domain in $\mathbb{R}^r$ with smooth boundary $\partial \Omega$, and
consider (\ref{rdnet})
where $x\in \mathbb{R}^n$, $f(\cdot)$ is a continuously differentiable vector field, and $\nabla^2x:=[\nabla^2x_1\cdots \nabla^2x_n]^T$ is the vector Laplacian.  In a typical reaction-diffusion system, $D\in \mathbb{R}^{n\times n}$ is a diagonal matrix of diffusion coefficients $d_i$ for each species $i=1,\cdots,n$; however, in the derivations below, we take $D$ to be an arbitrary real matrix for further generality.  We assume Neumann boundary conditions:
\begin{equation}
\nabla x_i(\xi) \cdot \hat{n}(\xi)=0 \quad \forall \xi \in \partial \Omega, \quad i=1,\cdots,n \label{bc}
\end{equation}
where $\xi$ represents the spatial variable and $\hat{n}$ is a vector normal to the boundary $\partial \Omega$. Well-posedness of (\ref{rdnet})-(\ref{bc}) is not emphasized in this paper; the reader may refer to \cite[Chapter 7.3]{smith} for conditions that guarantee existence of {\it classical} solutions to reaction-diffusion PDEs.

To establish a condition under which solutions $x(t,\xi)$ exhibit uniform behavior over the spatial domain $\Omega$, we
denote by:
\begin{equation}\label{pidef}
\pi\{v\}:=v-\bar{v}
\end{equation}
the deviation of a function $v=v(\xi)$ from its average:
\begin{equation}\label{avg}
\bar{v}:=\frac{1}{|\Omega|}\int_\Omega v(\xi)d\xi.
\end{equation}
In the derivations below, we also use the
$L_2(\Omega)$ inner product:
\begin{equation}
\langle u,v \rangle_{L_2(\Omega)}:=\int_\Omega u^T(\xi)v(\xi)d\xi
\end{equation}
and norm:
\begin{equation}
\|v\|_{L_2(\Omega)}:= \sqrt{\langle v,v\rangle_{L_2(\Omega)}}.
\end{equation}
We let $0=\lambda_1\le \lambda_2 \le \cdots \le \lambda_k\le \cdots$ denote the eigenvalues of the operator $L=-\nabla^2$ on $\Omega$ with Neumann boundary condition:
\begin{equation}\label{eigs}
L \phi_k(\xi)=\lambda_k \phi_k(\xi), \quad \nabla \phi_k(\xi) \cdot \hat{n}(\xi)=0  \ \, \forall \xi \in \partial \Omega,
\end{equation}
and make use of the second smallest eigenvalue, $\lambda_2$, in our main result:

\begin{theorem}\label{network}
Consider the reaction-diffusion system (\ref{rdnet})-(\ref{bc}) and let $\lambda_2$ be the second smallest eigenvalue of the operator $L=-\nabla^2$ on $\Omega$ with Neumann boundary condition as in (\ref{eigs}).  If there exists a convex set $\mathcal{X}\subseteq \mathbb{R}^n$, a matrix $P=P^T>0$, and a constant $\epsilon>0$ such that
\begin{eqnarray}\label{diagstab}
&&P\left(J(x)-\lambda_2D\right)+\left(J(x)-\lambda_2D\right)^TP\le -\epsilon I \qquad \forall x\in \mathcal{X}\\
&& PD+D^TP\ge 0, \label{diag}
\end{eqnarray}
then, for every classical solution $x(t,\xi):[0,\infty)\times \Omega \rightarrow \mathcal{X}$,
\begin{equation}
\|\pi\{x(t,\xi)\}\|_{L_2(\Omega)}\rightarrow 0
\end{equation}
 exponentially as $t\rightarrow \infty$. \hfill $\Box$
\end{theorem}

The second Neumann eigenvalue $\lambda_2$ is a measure of the well-connectedness of the spatial domain. Indeed,  of all sets of given volume, $\lambda_2$ is maximized by the ball \cite{henrot}. In situations where $\lambda_2$ is not easily calculable for the given domain $\Omega$, Theorem \ref{network} can be applied with a lower bound on $\lambda_2$ at the cost of making (\ref{diagstab})-(\ref{diag}) more restrictive.  A commonly used lower bound on $\lambda_2$ was derived for the Laplacian operator by Cheeger \cite{Che70}, and extended in \cite{chung} to Laplacian matrices of graphs.

Othmer \cite{othmer77}, followed by other papers \cite{AshOth78,ConHofSmo78,JonSle83},  studied the reaction-diffusion system (\ref{rdnet})-(\ref{bc}) with $D={\rm diag}\{d_1,\cdots,d_n\}$,
and proved uniform behavior of the solutions under the condition:
\begin{equation}\label{othm}
\sup_{x\in \mathcal{X}}\left\|J(x)\right\| < \lambda_2\min_i\{d_i\}.
\end{equation}
Note that (\ref{othm})  implies (\ref{diagstab}) with $P=I$, which means that Theorem \ref{network} incorporates Othmer's condition (\ref{othm}) as a special case.  Assumption (\ref{diagstab}) of Theorem \ref{network} is far less restrictive than (\ref{othm}),
 and is applicable to numerous practically important systems which do not satisfy global Lipschitz bounds.

\begin{example}\rm
As an illustration of Theorem \ref{network}, consider the Fitzhugh-Nagumo model of neuron excitation and oscillations (see {\it e.g.} \cite{keshet05}), augmented here with diffusion terms:
\begin{eqnarray}
\frac{\partial x_1}{\partial t}&=&c\left(x_1-\frac{1}{3}x_1^3+x_2\right)+d_1\nabla^2x_1\\
\frac{\partial x_2}{\partial t}&=&\frac{1}{c}\left(-x_1-bx_2+a\right)+d_2\nabla^2x_2, \quad c,b,d_1,d_2>0.
\end{eqnarray}
The Jacobian matrix:
\begin{equation}
J(x)=\left[\begin{array}{cc} c(1-x_1^2) & c \\ -\frac{1}{c} & -\frac{b}{c}
\end{array} \right]
\end{equation}
does not satisfy a norm bound; however, with $\lambda_2d_1>c$, conditions (\ref{diagstab})-(\ref{diag}) hold with:
\begin{equation}
P=\left[\begin{array}{cc} \frac{1}{c} & 0 \\ 0 & {c}
\end{array} \right].
\end{equation}
\mbox{}\hfill $\Box$
\end{example}

To see the implications of Theorem \ref{network} for a linear reaction-diffusion system, we recall that the solutions of (\ref{rdnet}) with $f(x)=Ax$ can be expanded  as:
\begin{equation}\label{expand}
x(t,\xi)=\sum_{k=1}^\infty \sigma_k(t)\phi_k(\xi)
\end{equation}
where, due to the orthogonality of the eigenfunctions $\phi_k(\xi)$ in (\ref{eigs}), $\sigma_k(t)\in \mathbb{R}^n$ obey the decoupled ODEs:
\begin{equation}
\dot{\sigma_k}=(A-\lambda_k D)\sigma_k.
\end{equation}
Because the eigenfunction $\phi_1$ corresponding to $\lambda_1=0$ is constant, the $k=1$ term  in (\ref{expand}) constitutes the average $\bar{x}$, governed by $\dot{\bar{x}}=A\bar{x}$, and the decay of the remaining terms is guaranteed when the matrices $A-\lambda_kD$, $k=2,3,\cdots,$ are Hurwitz. Conditions (\ref{diagstab})-(\ref{diag}) with $J(x)=A$ in Theorem \ref{network}  imply the existence of a common Lyapunov function for these matrices, thus ensuring that they are indeed Hurwitz.
\bigskip

\noindent
{\bf Proof of Theorem \ref{network}:}
We denote
\begin{equation}\label{dev}
\tilde{x}:=\pi\{x\}
\end{equation}
where $\pi\{\cdot\}$ is as defined in (\ref{pidef}), and note that $\tilde{x}$ satisfies:
\begin{equation}\label{rdnet2}
\frac{\partial \tilde{x}}{\partial t}=\pi\{f(x)\}+D\nabla^2x,
\end{equation}
where we have substituted $\pi\{\nabla^2x_i\}=\nabla^2x_i$ because $\int_\Omega \nabla^2x_i \,d\xi= \int_{\partial \Omega} \nabla x_i \cdot \hat{n}\,dS=0$
from the Divergence Theorem and the boundary condition (\ref{bc}).
We then select the functional:
\begin{equation}\label{lyap}
V(\tilde{x})=\frac{1}{2}\langle \tilde{x},P\tilde{x}\rangle_{L_2(\Omega)},
\end{equation}
where $P$ is as in (\ref{diagstab})-(\ref{diag}), and obtain:
\begin{equation}\label{claims2}
\dot{V} \le \langle \tilde{x}, P\pi\{f(x)\}\rangle_{L_2(\Omega)}+\langle \tilde{x}, PD\nabla^2x\rangle_{L_2(\Omega)}.
\end{equation}
We note from (\ref{diag}) that there exists a matrix $Q$ such that $Q^TQ=\frac{1}{2}(PD+D^TP)$.  This means that:
\begin{equation}\label{rmat}
\langle \tilde{x}, PD\nabla^2x\rangle_{L_2(\Omega)}=\langle Q\tilde{x}, Q\nabla^2\tilde{x}\rangle_{L_2(\Omega)}
=\langle y, \nabla^2y \rangle_{L_2(\Omega)},
\end{equation}
where $y:=Q\tilde{x}$.
Integrating both sides of  the identity
\begin{equation}
\nabla\cdot ({y}_i\nabla{y}_i)=|\nabla y_i|^2+y_i\nabla^2{y}_i
\end{equation}
over $\Omega$ and noting that the left-hand side vanishes due to the Divergence Theorem and the boundary condition (\ref{bc}), we obtain:
\begin{equation}\label{almost}
\int_\Omega y_i\nabla^2{y}_i d\xi=-\int_\Omega |\nabla{y}_i|^2 d\xi.
\end{equation}
Moreover, because $\int_\Omega y d\xi=Q\int_\Omega\tilde{x}d\xi=0$, it follows from the the Poincar\'{e} Inequality \cite[Equation (1.37)]{henrot}
that:
\begin{equation}\label{there}
  \int_\Omega |\nabla{y}_i|^2 d\xi \ge \lambda_2\int_\Omega {y}_i^2 d\xi
\end{equation}
and, thus, (\ref{almost}) and (\ref{there}) imply:
\begin{equation}\label{main1}
\langle y_i,\nabla^2 y_i \rangle_{L_2(\Omega)} \le -\lambda_2 \|y_i\|^2_{L_2(\Omega)}.
\end{equation}
We substitute the inequality (\ref{main1}) in (\ref{rmat}), substitute back $y=Q\tilde{x}$ and $Q^TQ=\frac{1}{2}(PD+D^TP)$, and obtain:
\begin{equation}\label{rmat2}
\langle \tilde{x}, PD\nabla^2x\rangle_{L_2(\Omega)}=\langle y, \nabla^2y \rangle_{L_2(\Omega)}\le -\lambda_2\langle y, y\rangle_{L_2(\Omega)}=-\lambda_2\langle \tilde{x}, PD \tilde{x}\rangle_{L_2(\Omega)}.
\end{equation}
Substitution of (\ref{rmat2}) in (\ref{claims2}) then gives:
\begin{equation}\label{claims3}
\dot{V} \le \langle \tilde{x}, P\pi\{f(x)\}\rangle_{L_2(\Omega)}-\lambda_2\langle \tilde{x}, PD \tilde{x}\rangle_{L_2(\Omega)}.
\end{equation}

Next, we rewrite the first term on the right-hand side of (\ref{claims3}) as:
\begin{eqnarray}\label{fir}
\langle \tilde{x}, P\pi\{f(x)\}\rangle_{L_2(\Omega)}&=&\int_\Omega \tilde{x}^TP\left(f(x)-\frac{1}{|\Omega|}\int_\Omega f(x) d\xi\right)d\xi \\ &=& \int_\Omega \tilde{x}^TP\left(f(x)-f(\bar{x})\right)d\xi+\int_\Omega \tilde{x}^TP\left(f(\bar{x})-\frac{1}{|\Omega|}\int_\Omega f(x) d\xi\right)d\xi \label{sec}
\\ &=& \int_\Omega \tilde{x}^TP\left(f(x)-f(\bar{x})\right)d\xi \label{thi}
\end{eqnarray}
where, to obtain (\ref{sec}), we added and subtracted $f(\bar{x})$ in (\ref{fir}).  To obtain (\ref{thi}),  we noted that the second integral in (\ref{sec}) is zero because the factor
\begin{equation}
\left(f(\bar{x})-\frac{1}{|\Omega|}\int_\Omega f(x) d\xi\right)
\end{equation}
does not depend on $\xi$, and because $\int_\Omega \tilde{x}d\xi=0$.  Substitution of (\ref{thi}) in (\ref{claims3}) then results in:
\begin{equation}\label{claims4}
\dot{V} \le \langle \tilde{x}, P(f(x)-f(\bar{x})-\lambda_2D\tilde{x})\rangle_{L_2(\Omega)}.
\end{equation}

Finally, we use the Mean-Value Theorem \cite{ortegarh} and write:
\begin{equation}\label{MV}
f(x)-f(\bar{x})=\int_0^1 J(\bar{x}+s(x-\bar{x}))(x-\bar{x})ds.
\end{equation}
Substituting in (\ref{claims4}) and using (\ref{diagstab}), we obtain:
\begin{equation}\label{last}
\dot{V}\le \int_0^1\int_\Omega \tilde{x}^TP\left(J(\bar{x}+s\tilde{x})-\lambda_2D\right)\,\tilde{x}\,d\xi \,ds \le \int_0^1 \int_\Omega -\frac{\epsilon}{2} \tilde{x}^T\tilde{x}\,d\xi \, ds\le -\frac{\epsilon}{\lambda_{\max} (P)}V.
\end{equation}
Inequality (\ref{last}) proves exponential decay of the functional $V(\tilde{x})$ defined in (\ref{lyap}), from which the conclusion of the theorem follows. \hfill $\Box$

\section{Constant Matrix Parameterizations of the Jacobian}\label{LMIsec}
We now present a procedure to verify (\ref{diagstab}) by bounding the Jacobian $J(x)$ within a set which is parameterized by constant matrices.  Examples of such parameterizations include the convex hull:
\begin{equation}\label{conv}
conv\{Z_1,\cdots,Z_q\}=\left\{\, \theta_1Z_1+\cdots+\theta_qZ_q \ | \ \theta_1+\cdots+\theta_q=1, \ \theta_i\ge 0 \  \,i=1,\cdots,q\right\},
\end{equation}
and the conic hull:
\begin{equation}\label{cone}
cone\{S_1,\cdots,S_m\}=\left\{\, \omega_1S_1+\cdots+\omega_mS_m \ | \ \omega_i\ge 0 \  \,i=1,\cdots,m\right\}.
\end{equation}
When $J(x)$ belongs to the sum of these two sets, (\ref{diagstab}) can be replaced with the constant matrix inequalities (\ref{hull1})-(\ref{hull2}) below:
\begin{theorem}\label{LMIlem}
If there exist constant matrices $Z_1,\cdots,Z_q$ and $S_1,\cdots,S_m$ such that
\begin{equation}\label{membership}
J(x) \in conv\{Z_1,\cdots,Z_q\}+cone\{S_1,\cdots,S_m\} \quad \forall x\in \mathcal{X},
\end{equation}
then a matrix $P=P^T$ satisfying:
\begin{eqnarray}\label{hull1}
P(Z_k-\lambda_2D)+(Z_k-\lambda_2D)^TP<0, &&\quad k=1,\cdots,q \\
PS_k+S_k^TP\le 0, &&\quad k=1,\cdots,m \label{hull2}
\end{eqnarray}
also satisfies (\ref{diagstab}) for some $\epsilon>0$.
If the image of $\mathcal{X}$ under $J(\cdot)$ is surjective onto $$conv\{Z_1,\cdots,Z_q\}+cone\{S_1,\cdots,S_m\},$$ then the converse is also true; that is, (\ref{diagstab}) with $\epsilon>0$ implies (\ref{hull1})-(\ref{hull2}). \hfill $\Box$
\end{theorem}

The proof is routine and is given in the Appendix.  Theorem \ref{LMIlem} is useful because the inequalities (\ref{diag}), (\ref{hull1}) and (\ref{hull2}) are linear in the variables $P=P^T>0$ and $\epsilon>0$ and, thus, the conditions of Theorem \ref{network}
 can be checked with efficient numerical tools available for {\it linear matrix inequalities} \cite{boyd}.  Analytical conditions for the existence of common quadratic Lyapunov functions are also available for several classes of matrices \cite{LibMor99}.

In various examples of reaction networks, $J(x)$ belongs to a convex set of the form:
\begin{equation}\label{boxset}
box\{A_0,A_1, \cdots,A_\ell\}=\{A_0+\gamma_1A_1+\cdots+\gamma_\ell A_\ell\ | \ 0\le \gamma_i \le 1, \ i=1,\cdots,\ell\},
\end{equation}
where $A_1,\cdots,A_\ell$ are rank-one matrices.  Although Theorem \ref{LMIlem} is applicable to the matrices $Z_1,\cdots,Z_q$ corresponding to the vertices of the set (\ref{boxset}), this application involves $q=2^\ell$ vertices and may become intractable for large $\ell$. Theorem \ref{dyad} below, proven in the Appendix, gives an alternative test  that uses only the matrices $A_0,\cdots,A_\ell$ for verifying (\ref{diagstab}):

\begin{theorem}\label{dyad}
Suppose
\begin{equation}\label{alter}
J(x)\in box\{A_0,A_1, \cdots,A_\ell\},
\end{equation}
where $A_1, \cdots, A_\ell$ are rank-one matrices and, thus, can be decomposed as:
\begin{equation}\label{dya}
A_i=B_iC_i^T \quad i=1,\cdots,\ell
\end{equation}
with appropriately selected column vectors $B_i,C_i \in \mathbb{R}^{n}$.
If there exists a matrix $\mathcal P=\mathcal P^T>0$ of the form:
\begin{equation}\label{blockdia}
\mathcal{P}=\left[\begin{array}{cccc}P &  & &  \\  & q_1 &  & \\ & &  \ddots & \\ & & & q_\ell \end{array}\right],\quad P\in \mathbb{R}^{n\times n}, \ q_i\in \mathbb{R}, \ i=1,\cdots,\ell,
\end{equation}
 such that:
\begin{equation}\label{composite}
{\mathcal P}\left[\begin{array}{cc} A_0-\lambda_2 D &  B \\ C^T & -I\end{array}\right]+\left[\begin{array}{cc} A_0-\lambda_2 D &  B \\ C^T & -I\end{array}\right]^T{\mathcal P}<0,
\end{equation}
where $B:=[B_1\cdots B_\ell]$ and $C:=[C_1\cdots C_\ell]$, then the upper-left block $P=P^T>0$ satisfies (\ref{diagstab}) for some $\epsilon>0$. If, in addition, $\ell=1$ and the image of $\mathcal{X}$ under $J(\cdot)$ is surjective  onto $box\{A_0,A_1\}$, then  the converse is also true; that is, if (\ref{diagstab}) holds with  a matrix $P=P^T>0$ and a constant $\epsilon>0$, then there exists $q_1>0$ such that $\mathcal P=\mathcal P^T>0$ in (\ref{blockdia}) satisfies (\ref{composite}).\hfill $\Box$
\end{theorem}

\begin{remark}\label{diagonal}\rm  In applications, it may be preferable to search for a fully diagonal matrix $\mathcal{P}$ satisfying (\ref{composite}), instead of a block-diagonal $\mathcal{P}$ as in (\ref{blockdia}).  Although this  restriction may add conservatism, it has the following advantages:

i) Condition (\ref{diag}) in Theorem \ref{network} is satisfied for all diagonal and nonnegative $D$, and need not be checked separately when $D$ has this form.

ii) If the set $box\{A_0,A_1, \cdots,A_\ell\}$ in Theorem \ref{dyad} is augmented with
$cone\{S_1,\cdots,S_m\},$  then  the upper left $n\times n$ component $P$ of the matrix $\mathcal{P}$ in (\ref{blockdia}) must satisfy (\ref{hull2}) in addition to (\ref{composite}).  However, in the special case where $S_k$, $k=1,\cdots,m,$ are nonpositive diagonal matrices,  (\ref{hull2}) holds for every diagonal $P>0$ and, thus, it is sufficient to check (\ref{composite}) with a diagonal $\mathcal{P}>0$.

iii) Likewise, if some of the matrices $A_i$ in (\ref{alter}) are diagonal and nonpositive, the corresponding columns $B_i$ and $C_i$ can be omitted in constructing the matrix
\begin{equation}\label{redf}
\mathcal{A}=\left[\begin{array}{cc} A_0-\lambda_2 D &  B \\ C^T & -I\end{array}\right],
\end{equation}
thus reducing the dimension of the problem (\ref{composite}).

iv) In several practically important examples, analytical tests are applicable to check the existence of a diagonal solution to the Lyapunov inequality (\ref{composite}).  Matrices $\mathcal{A}$ for which a diagonal $\mathcal{P}>0$ satisfying
\begin{equation}
\mathcal{PA}+\mathcal{A}^T\mathcal{P}<0
\end{equation}
exists are termed {\it diagonally stable} \cite{bhaya}, and have been fully characterized in dimension three \cite{Cro78} and dimension four \cite{Red85}. For higher dimensional matrices, diagonal stability tests have been derived by exploiting special sparse structures, such as a {\it cyclic} structure and its variants studied in \cite{ArcSon06,ArcSon08}.
If the matrix (\ref{redf})
conforms to one of these structures, the existence of a diagonal solution $\mathcal{P}$ to (\ref{composite}) can be checked with simple algebraic conditions.
  Analytical conditions are indeed important in applications, because they reveal which system properties and which parameters are critical for Theorem \ref{network} to hold.
\hfill $\Box$
\end{remark}

\begin{example}\label{Goodwinex}\rm In \cite{JacMon61}, Jacob and Monod gave a molecular description of how certain metabolites regulate their production by repressing enzymes necessary for their synthesis.  Following this description, Goodwin \cite{Goodwin65} proposed a differential equation model and studied its oscillatory behavior.
 A variant of Goodwin's model \cite{Thr91},  augmented here with diffusion terms, is:
\begin{eqnarray}\nonumber
\frac{\partial x_1}{\partial t}&=& -a_1x_1+\frac{V_1}{K_1+x_3}+d_1\nabla^2x_1\\
\frac{\partial x_2}{\partial t}&=&-a_2x_2+b_1x_1+d_2\nabla^2x_2 \label{good}\\
\frac{\partial x_3}{\partial t}&=&-\frac{V_3x_3}{K_3+x_3}+b_2x_2+d_3\nabla^2x_3, \nonumber
\end{eqnarray}
where all parameters are positive and $x_1$, $x_2$, $x_3$ denote, respectively, the concentrations of the messenger RNA, enzyme and product.

To inspect condition (\ref{diagstab}) of Theorem \ref{network} on the set $\mathcal{X}=\mathbb{R}^3_{\ge 0}$,
we study the Jacobian matrix:
\begin{equation}
J(x)=\left[\begin{array}{ccc} -a_1 & 0 & -b_3(x_3)\\ b_1 & -a_2 & 0 \\
0 & b_2 & -a_3(x_3) \end{array}\right],
\end{equation}
and note that
\begin{equation}
a_3(x_3):=\frac{V_3K_3}{(K_3+x_3)^2} \qquad b_3(x_3):=\frac{V_1}{(K_1+x_3)^2}
\end{equation}
lie in the bounded intervals $[0,\frac{V_3}{K_3}]$ and $[0,\frac{V_1}{K_1^2}]$, respectively.
This means that
\begin{equation}\label{boxsetex}
J(x)\in box\{A_0,A_1,A_2\},
\end{equation}
where
\begin{equation}
A_0=\left[\begin{array}{ccc} -a_1 & 0 & 0\\ b_1 & -a_2 & 0 \\
0 & b_2 & 0 \end{array}\right] \quad  A_1=\left[\begin{array}{ccc} 0 & 0 & -\frac{V_1}{K_1^2}\\ 0 & 0 & 0 \\
0 & 0 & 0 \end{array}\right]\quad  A_2=\left[\begin{array}{ccc} 0 & 0 & 0\\ 0 & 0 & 0 \\
0 & 0 & -\frac{V_3}{K_3} \end{array}\right].
\end{equation}
We decompose $A_1$ as $A_1=B_1C_1^T$ with $B_1=[-\frac{V_1}{K_1^2} \ 0 \ 0]^T$ and $C_1=[0 \ 0 \ 1]^T$, and construct the matrix:
\begin{equation}\label{wa}
\left[\begin{array}{cc} A_0-\lambda_2 D &  B_1 \\ C_1^T & -I\end{array}\right]=\left[\begin{array}{cccc} -(a_1+\lambda_2d_1) & 0 & 0 & -\frac{V_1}{K_1^2}\\ b_1 & -(a_2+\lambda_2d_2) & 0 & 0 \\
0 & b_2 & -\lambda_2d_3 & 0 \\ 0 & 0 & 1 & -1\end{array}\right],
\end{equation}
where we have omitted $A_2$ in view of item (iii) in Remark \ref{diagonal}. Because the matrix (\ref{wa}) has a cyclic form, the {\it secant criterion} derived in \cite{ArcSon06} is applicable, and states that diagonal stability of (\ref{wa}) is equivalent to the condition:
\begin{equation}\label{seccon}
\frac{b_1b_2V_1}{K_1^2(a_1+\lambda_2d_1)(a_2+\lambda_2d_2)\lambda_2d_3}<\sec(\pi/4)^4=4.
\end{equation}
We thus conclude from Theorem \ref{dyad} and Remark \ref{diagonal} that, if the parameters of the model (\ref{good}) are such that (\ref{seccon}) holds with $\lambda_2$ calculated from the domain $\Omega$, then Theorem \ref{network} guarantees spatial uniformity of the solutions.  

Note that, in this example, $\|J(x)\|$ is bounded and, hence, condition (\ref{othm}) of \cite{othmer77} is applicable.   With the following set of parameters from \cite{Thr91}:
\begin{equation}\label{params}
a_1=a_2=b_1=b_2=0.01, \ V_1=9, \ V_3=K_1=K_3=1,
\end{equation}
(\ref{othm}) stipulates:
\begin{equation}\label{othmex}
\lambda_2\min_i{d_i}>\sup_{x\in \mathcal{X}}\left\|J(x)\right\|=9.0554,
\end{equation}
 where the $\sup$ is achieved when $x_3=0$. To compare this condition to (\ref{seccon}), we note that, for the same parameter values,
 \begin{equation}\label{upbo}
 \frac{b_1b_2V_1}{K_1^2(a_1+\lambda_2d_1)(a_2+\lambda_2d_2)\lambda_2d_3}<\frac{9\cdot10^{-4}}
 {(0.01+\lambda_2\min_i{d_i})^2\lambda_2\min_i{d_i}},
 \end{equation}
 which implies that (\ref{seccon}) holds if the upper-bound in (\ref{upbo}) is less than $4$; that is, if:
\begin{equation}\label{myex}
\lambda_2\min_i{d_i}>{0.05435}.
\end{equation}
The estimate (\ref{myex}) is obtained using the upper-bound (\ref{upbo}), which is achieved only when the diffusion coefficients are identical.  For nonidentical diffusion coefficients, condition (\ref{seccon}) leads to even more dramatic improvements over the conservative estimate (\ref{othmex}). In this example, a direct application of Theorem \ref{LMIlem} to the vertices of the set (\ref{boxsetex}), without insisting on a diagonal solution $P$, gave an insignificant improvement over (\ref{myex}):  With $d_1=d_2=d_3=d$, we numerically obtained the bound $\lambda_2d>0.05425$.
\end{example}

\begin{example}\label{Goldbeterex} \rm As a further illustration of Theorem \ref{dyad}, we consider a model of {\it Drosphila} circadian rhythms, proposed in \cite{gold95} and further studied in \cite{AngSon08} for its dynamical behavior.  When augmented with diffusion terms, this model is of the form:
\begin{eqnarray}
\frac{\partial M}{\partial t}&=& \frac{v_sK_I^n}{K_I^n+P_N^n}-\frac{v_mM}{k_m+M}+d_M\nabla^2M\\
\frac{\partial {P}_0}{\partial t}&=&k_sM-\frac{V_1P_0}{K_1+P_0}+\frac{V_2P_1}{K_2+P_1}+d_{P_0}\nabla^2P_0\\
\frac{\partial {P}_1}{\partial t}&=&\frac{V_1P_0}{K_1+P_0}-\frac{V_2P_1}{K_2+P_1}-\frac{V_3P_1}{K_3+P_1}+\frac{V_4P_2}{K_4+P_2}+d_{P_1}\nabla^2P_1\\
\frac{\partial {P}_2}{\partial t}&=&\frac{V_3P_1}{K_3+P_1}-\frac{V_4P_2}{K_4+P_2}-k_1P_2+k_2P_N-\frac{v_dP_2}{k_d+P_2}+d_{P_2}\nabla^2P_2\\
\frac{\partial {P}_N}{\partial t}&=&k_1P_2-k_2P_N+d_{P_N}\nabla^2P_N,
\end{eqnarray}
where $P_i$ represents the concentration of the PER protein, with the indices $i=0,1,2$ denoting the degree of phosphorylation.  Likewise, $P_N$ is the concentration of PER in the nucleus and $M$ is the concentration of the messenger RNA.

We obtain the Jacobian matrix
\begin{equation}\label{goldjaco}
J=\left[ \begin{array}{ccccc} -\phi_6(M)& 0 & 0 & 0 & -\phi_5(P_N) \\
k_s & -\phi_1(P_0) & \phi_2(P_1) & 0 & 0 \\
0 & \phi_1(P_0) &  -\phi_2(P_1)-\phi_3(P_1) & \phi_4(P_2) & 0\\
0 & 0 & \phi_3(P_1) & -k_1-\phi_4(P_2)-\phi_7(P_2) & k_2
\\ 0 & 0 & 0 & k_1 & -k_2
\end{array}\right]
\end{equation}
where
\begin{equation}
\phi_i(x):=\frac{K_iV_i}{(K_i+x)^2} \quad i=1,2,3,4, \quad \phi_5(x):=\frac{nv_sK_I^nx^{n-1}}{(K_I^n+x^n)^2}, \quad \phi_6(x):=\frac{v_mk_m}{(k_m+x)^2},\quad \phi_7(x):=\frac{v_dk_d}{(k_d+x)^2},
\end{equation}
and note that, for $x\ge 0$, these functions lie in the bounded intervals $[0,\bar{\phi}_i]$, with:
\begin{equation}
\bar{\phi}_i=\frac{V_i}{K_i}\quad i=1,2,3,4, \quad \bar{\phi}_6=\frac{v_m}{k_m}, \quad \bar{\phi}_7=\frac{v_d}{k_d}, \quad \bar{\phi}_5=\left\{ \begin{array}{ll} \frac{nv_s}{K_I^{n}} & \mbox {if} \ n=1 \\ \frac{(n+1)^2v_s}{4nK_I^{n}}\left(K_I^n\frac{n-1}{n+1}\right)^{\frac{n-1}{n}} & \mbox {if} \ n>1.\end{array} \right.
\end{equation}
Thus, for all $(M,P_0,P_1,P_2,P_N)\in \mathbb{R}^5_{\ge 0},$ $J\in box\{A_0,A_1,\cdots,A_7\}$, where
\begin{equation}
A_0=\left[ \begin{array}{ccccc} 0 & 0 & 0 & 0 & 0 \\
k_s & 0 & 0 & 0 & 0 \\
0 & 0 &  0 & 0 & 0\\
0 & 0 & 0 & -k_1 & k_2
\\ 0 & 0 & 0 & k_1 & -k_2
\end{array}\right]\ A_1=\left[ \begin{array}{ccccc} 0 & 0 & 0 & 0 & 0 \\
0 & -\frac{V_1}{K_1} & 0 & 0 & 0 \\
0 &  \frac{V_1}{K_1} & 0  & 0 & 0\\
0  & 0 & 0 & 0 & 0
\\ 0 & 0 & 0 & 0 & 0
\end{array}\right]\ A_2=\left[ \begin{array}{ccccc} 0 & 0 & 0 & 0 & 0 \\
0 & 0 & \frac{V_2}{K_2} & 0 & 0 \\
0 &  0 & -\frac{V_2}{K_2}  & 0 & 0\\
0  & 0 & 0 & 0 & 0
\\ 0 & 0 & 0 & 0 & 0
\end{array}\right]
\end{equation}
\begin{equation}
A_3=\left[ \begin{array}{ccccc} 0 & 0 & 0 & 0 & 0 \\
0 & 0 & 0 & 0 & 0 \\
0 &  0 & -\frac{V_3}{K_3}  & 0 & 0\\
0  & 0 & \frac{V_3}{K_3} & 0 & 0
\\ 0 & 0 & 0 & 0 & 0
\end{array}\right]\ A_4=\left[ \begin{array}{ccccc} 0 & 0 & 0 & 0 & 0 \\
0 & 0 & 0 & 0 & 0 \\
0 &  0 & 0 & \frac{V_4}{K_4}  & 0 \\
0  & 0 & 0 & -\frac{V_4}{K_4} & 0
\\ 0 & 0 & 0 & 0 & 0
\end{array}\right]\ A_5=\left[ \begin{array}{ccccc} 0 & 0 & 0 & 0 & -\bar{\phi}_5 \\
0 & 0 & 0 & 0 & 0 \\
0 &  0 & 0 & 0  & 0 \\
0  & 0 & 0 & 0 & 0
\\ 0 & 0 & 0 & 0 & 0
\end{array}\right]
\end{equation}
\begin{equation}
A_6=\left[ \begin{array}{ccccc} -\frac{v_m}{k_m} & 0 & 0 & 0 & 0 \\
0 & 0 & 0 & 0 & 0 \\
0 &  0 & 0  & 0 & 0\\
0  & 0 & 0 & 0 & 0
\\ 0 & 0 & 0 & 0 & 0
\end{array}\right]\ A_7=\left[ \begin{array}{ccccc} 0 & 0 & 0 & 0 & 0 \\
0 & 0 & 0 & 0 & 0 \\
0 &  0 & 0 & 0  & 0 \\
0  & 0 & 0 & -\frac{v_d}{k_d} & 0
\\ 0 & 0 & 0 & 0 & 0
\end{array}\right].
\end{equation}
Using the following parameter values from \cite{gold95}:
\begin{eqnarray}\nonumber
&& n=4,\ v_s=0.76,\ K_I=1,\ k_s=0.38,\ k_1=1.9,\ k_2=1.3,\ V_1=3.2,\ V_2=1.58,\ V_3=5,\ V_4=2.5,\\
&& K_1=2,\ K_2=2,\ K_3=2,\ K_4=2,\ v_m=0.65,\ v_d=0.95,\ k_d= 0.2,\ k_m=0.5,
\end{eqnarray}
and assuming identical diffusion coefficients, denoted by $d$, we applied the procedure outlined in Theorem \ref{dyad} and numerically determined bounds for $\lambda_2d$ using the MATLAB software CVX \cite{GraBoy08}.  The linear matrix inequality (\ref{blockdia})-(\ref{composite}) was feasible with a fully populated matrix $P$ when $\lambda_2d\ge 0.4590$, and with a diagonal matrix $P$ when $\lambda_2d\ge 0.5393$.

Note that, in our parameterization, we took advantage of the repetition of the nonlinearities in $\phi_1,\cdots \phi_4$ in (\ref{goldjaco}), and employed the matrices $A_1,\cdots,A_4$, each representing two occurrences of the same nonlinearity.  The alternative approach of overparameterizing with one matrix for each occurrence would lead to conservative results.  Indeed, a repetition of the numerical experiment described above with $A_1,\cdots,A_4$ split into two matrices each, gave the conservative feasibility region $\lambda_2d\ge 1.7892$ with a fully populated $P$ (compare to $\lambda_2d\ge 0.4590$ above), and $\lambda_2d\ge 1.7943$ with a diagonal $P$ (compare to $\lambda_2d\ge 0.5393$ above).
\end{example}

\section{Synchronization in a Network of ODEs with Diffusion-Like Coupling}\label{odesec}
We now derive an analogous result for a network of identical ODE models that are interconnected according to an undirected graph:
\begin{equation}\label{initial}
\dot{x}^k=f(x^k)+D\sum_{j\in \mathcal{N}_k}(x^j-x^k) \quad k=1,\cdots,N,
\end{equation}
where $x^k\in \mathbb{R}^n$, $\mathcal{N}_k\subseteq \{1,\cdots,N\}$ denotes the set of nodes adjacent to node $k$ in the graph, and  $D$ is allowed to be an arbitrary real matrix as in Section \ref{pdesec}. Denoting by $X$ the concatenated vector:
\begin{equation}
X=[{x^1}^T \cdots {x^N}^T]^T,
\end{equation}
and by $L=(l_{i,j})\in\mathbb{R}^{N\times N}$ the graph Laplacian matrix \cite{godsil}:
\begin{equation}
l_{i,j}=\left\{ \begin{array}{l} \text{number of nodes adjacent to node $i$ if $i=j$} \\
\text{$-1$ if $i\neq j$ and $j\in \mathcal{N}_i$} \\ \text{$0$ otherwise,}\end{array} \right.
\end{equation}
we rewrite (\ref{initial}) in the compact form:
\begin{equation}\label{compa}
\dot{X}={F}(X)-(L\otimes D)X,
\end{equation}
where ``$\otimes$" represents the Kronecker product, and
\begin{equation}\label{Fnot}
F(X):=[f(x^1)^T \cdots f(x^N)^T]^T.
\end{equation}

We let $0=\lambda_1\le \lambda_2\le \cdots\le \lambda_N$ denote the eigenvalues of the Laplacian matrix, and show that the components $x^k(t)$ in (\ref{initial}) synchronize if $\lambda_2$ is such that (\ref{diagstab})-(\ref{diag}) hold as in Theorem \ref{network}:

\begin{theorem}\label{odethm}
Consider the interconnected system (\ref{compa})-(\ref{Fnot}), and suppose (\ref{diagstab})-(\ref{diag}) hold with a  matrix $P=P^T>0$ and a constant $\epsilon>0$ on a convex set $\mathcal{X}\subseteq \mathbb{R}^n$.  Then, every forward-complete solution $X(t)=[x^1(t)^T \cdots x^N(t)^T]^T$ that remains in $\mathcal{X}^N$ has the property that, for any pair $(k,j)\in \{1,\cdots,N\}\times \{1,\cdots,N\}$,
\begin{equation}\label{sync}
x^k(t)-x^j(t)\rightarrow 0
\end{equation}
exponentially as $t\rightarrow \infty$. \hfill $\Box$
\end{theorem}

\noindent
{\bf Proof of Theorem \ref{odethm}:} In this proof, we make repeated use of the property:
\begin{equation}\label{identity}
(A\otimes B)(C\otimes D)=(AC)\otimes (BD),
\end{equation}
which holds whenever the matrices are of compatible dimensions to form the indicated products.  We also recall that the Laplacian matrix $L$ satisfies:
\begin{equation}\label{null}
L1_N=0
\end{equation}
where $1_N$ denotes the $N\times 1$ vector of ones. Since $1_N$ is an eigenvector that corresponds to the eigenvalue $\lambda_1=0$, it follows that
\begin{equation}\label{luca}
y^TLy\ge \lambda_2 y^Ty \quad \forall y \perp 1_N.
\end{equation}
Likewise, denoting by $I_n$ the $n\times n$ identity matrix, we get the inequality:
\begin{equation}\label{akin}
y^T(L\otimes I_n) y\ge \lambda_2 y^Ty \quad \forall y \perp 1_N\otimes I_n,
\end{equation}
which is the discrete analog of the Poincar\'{e} Inequality (\ref{there}) used in the proof of Theorem \ref{network}.

Mimicking (\ref{avg}) and (\ref{dev}), we define:
\begin{equation}
\bar{x}:=\frac{1}{N}(x^1+\cdots+x^N)=\frac{1}{N}(1_N^T\otimes I_n)X, \quad \bar{X}:=1_N\otimes \bar{x},
\end{equation}
and
\begin{equation}
\tilde{x}^k:=x^k-\bar{x}, \quad \tilde{X}:=X-\bar{X}.
\end{equation}
It follows from this definition that $\sum_{k=1}^N\tilde{x}^k=0$ and, thus, for any matrix $M$ with $n$ rows,
\begin{equation}\label{M}
\tilde{X}^T(1_N\otimes M)=\sum_{k=1}^N\tilde{x}^{k}{}^TM=0.
\end{equation}
The dynamics of $\tilde{X}$ are given by:
\begin{eqnarray}\nonumber
\dot{\tilde{X}}&=&F(X)-\dot{\bar{X}}-(L\otimes D)X\\
&=&F(X)-\dot{\bar{X}}-(L\otimes D)\tilde{X},
\end{eqnarray}
where the second equation follows by substituting $X=\tilde{X}+1_N\otimes \bar{x}$ and by noting from (\ref{identity}) and (\ref{null}) that $(L\otimes D)(1_N\otimes \bar{x})$$=(L1_N)\otimes (D\bar{x})=0$.

We introduce the Lyapunov function $V=\frac{1}{2}\tilde{X}^T(I_N\otimes P)\tilde{X}$ and note that it satisfies:
\begin{eqnarray}\nonumber
\dot{V}&=&\tilde{X}^T(I_N\otimes P)(F(X)-\dot{\bar{X}})-\tilde{X}^T(I_N\otimes P)(L\otimes D)\tilde{X}\\
&=&\tilde{X}^T(I_N\otimes P)(F(X)-\dot{\bar{X}})-\tilde{X}^T(L\otimes (PD))\tilde{X}.\label{start}
\end{eqnarray}
Because $L$ is symmetric, the following identity holds:
\begin{equation}\label{clef}
(L\otimes (PD))+(L\otimes (PD))^T=L\otimes (PD+D^TP).
\end{equation}
As in the proof of Theorem \ref{network}, we define $Q$ such that $Q^TQ=\frac{1}{2}(PD+D^TP)$ and obtain:
\begin{equation}\label{bab}
\tilde{X}^T(L\otimes (PD))\tilde{X}=\tilde{X}^T(I_N\otimes Q^T)(L\otimes I_n)(I_N\otimes Q)\tilde{X}=y^T(L\otimes I_n)y,
\end{equation}
where
\begin{equation}\label{ydef}
y:=(I_N\otimes Q)\tilde{X}.
\end{equation}
It then follows from (\ref{identity}) and (\ref{M}) with $M=Q^T$ that
\begin{equation}
y^T(1_N\otimes I_n)=\tilde{X}^T(I_N\otimes Q^T)(1_N\otimes I_n)=\tilde{X}^T(1_N\otimes Q^T)=0,
\end{equation}
which means $y \perp 1_N\otimes I_n$ and, thus, the inequality (\ref{akin}) above is applicable. Using (\ref{akin}), (\ref{bab}) and (\ref{ydef}), we obtain:
\begin{equation}
\tilde{X}^T(L\otimes (PD))\tilde{X}=y^T(L\otimes I_n)y \ \ge \ \lambda_2 y^Ty=\lambda_2\tilde{X}^T(I_N\otimes (PD))\tilde{X}=\lambda_2\sum_{k=1}^N \tilde{x}^k {}^T\!\!PD\tilde{x}^k.
\end{equation}
Substituting this inequality back in (\ref{start}), we get:
\begin{equation}\label{f7}
\dot{V}\le \tilde{X}^T(I_N\otimes P)(F(X)-\dot{\bar{X}})-\lambda_2\sum_{k=1}^N \tilde{x}^k {}^T\!\!PD\tilde{x}^k.
\end{equation}

We next add and subtract $F(\bar{X})=1_N\otimes f(\bar{x})$, and rewrite (\ref{f7}) as:
\begin{eqnarray}\nonumber
\dot{V}&\le& \tilde{X}^T(I_N\otimes P)(F(X)-F(\bar{X}))+\tilde{X}^T(I_N\otimes P)(1_N\otimes(f(\bar{x})-\dot{\bar{x}}))-\lambda_2\sum_{k=1}^N \tilde{x}^k {}^T\!\!PD\tilde{x}^k\\ \nonumber
&=&\tilde{X}^T(I_N\otimes P)(F(X)-F(\bar{X}))+\tilde{X}^T(1_N\otimes P(f(\bar{x})-\dot{\bar{x}}))-\lambda_2\sum_{k=1}^N \tilde{x}^k {}^T\!\!PD\tilde{x}^k\\
&=&\tilde{X}^T(I_N\otimes P)(F(X)-F(\bar{X}))-\lambda_2\sum_{k=1}^N \tilde{x}^k {}^T\!\!PD\tilde{x}^k, \label{f8}
\end{eqnarray}
where the second equation follows from (\ref{identity}) and the third equation follows from (\ref{M}) with $M=P(f(\bar{x})-\dot{\bar{x}})$. Expanding the first term in (\ref{f8}) as a summation, we obtain:
\begin{equation}
\dot{V}\le \sum_{k=1}^N \tilde{x}^k {}^T\!\!P(f(x^k)-f(\bar{x}))-\lambda_2\sum_{k=1}^N \tilde{x}^k {}^T\!\!PD\tilde{x}^k.
\end{equation}
Finally, an application of the Mean-Value Theorem (\ref{MV}) yields:
\begin{equation}
\dot{V}\le \sum_{k=1}^N\int_0^1 \tilde{x}^k {}^T\!\!P(J(\bar{x}+s\tilde{x}_i)-\lambda_2D)\tilde{x}^k\,ds \le -\frac{\epsilon}{2} \tilde{X}^T\tilde{X}\le -\frac{\epsilon}{\lambda_{\max}(P)} V,
\end{equation}
which concludes the proof.  \hfill $\Box$
\begin{remark}\label{directed}\rm
In Theorem \ref{odethm}, we assumed an undirected graph to give an exact analog to the reaction-diffusion PDE result of Theorem \ref{network}.  However,  with the additional condition that the product $PD$ be symmetric, it is not difficult to extend Theorem \ref{odethm} to a directed graph, where $L$ is not symmetric and is restricted only by $L1_N=0$. In this extension, (\ref{clef}) must be replaced with:
\begin{equation}
(L\otimes (PD))+(L\otimes (PD))^T=(L+L^T)\otimes (PD)
\end{equation}
which holds because $PD$ is symmetric, and $\lambda_2$ must be redefined as the largest number such that (\ref{luca}) holds.  This definition of $\lambda_2$ was introduced in \cite{Wu05} as the ``algebraic connectivity" of a directed graph, and employed in \cite{ScaArcSon09} to obtain a synchronization result over directed and weighted graphs.\hfill $\Box$
\end{remark}

It is important to note that the Lyapunov inequalities (\ref{diagstab})-(\ref{diag}) used in Theorems \ref{network} and \ref{odethm} imply a {\it contraction} property \cite{LS98} for the family of vector fields $\{f(x)-\lambda Dx, \ \lambda\ge \lambda_2\}$.  Contraction properties, in various forms, have been employed in \cite{WuChu95,PogNij01,Wu05b,StaSep07,ScaArcSon09,RusBer09} to derive synchronization conditions for networks. We now make a connection between Theorem \ref{odethm} and one of the results in \cite{StaSep07}:

\begin{example}\label{rodolphe} \rm  Stan and Sepulchre \cite{StaSep07} studied the ODE models\footnote{We follow a slightly different notation than \cite{StaSep07} for consistency with Theorem \ref{odethm}.}:
\begin{eqnarray}\label{stan1}
\dot{x}^k&=&Ax^k+B\phi(y^k)+Bu^k\\
y^k&=&Cx^k,\label{stan2}
\end{eqnarray}
$x^k\in \mathbb{R}^n$, $u^k\in \mathbb{R}$, $y^k \in \mathbb{R}$, $k=1,\cdots,N$,
coupled by the feedback law
\begin{equation}
u=-Ly,\label{stan3}
\end{equation}
where $u:=[u^1\cdots u^N]^T$ and $y:=[y^1\cdots y^N]^T$, and pursued {\it incremental passivity} arguments to prove synchronization of the subsystems.

We now show that Theorem \ref{odethm}  is applicable to (\ref{stan1})-(\ref{stan3}) when the following hypotheses, adapted\footnote{Unlike \cite{StaSep07}, in H1, we assume that the nonlinearity $\phi(\cdot)$ is differentiable.  In H2, we strengthen the passivity assumption of \cite{StaSep07} to strict passivity so that (\ref{kyp}) holds with strict inequality and, thus, Theorem \ref{odethm} is directly applicable. With a slight modification of Theorem \ref{odethm}, it is indeed possible to remove the strictness condition and, instead, to assume observability of the pair $(A,C)$ as in \cite{StaSep07}. Finally, in H3, we remove the ``balanced graph" assumption ($1_N^TL=0$) employed in \cite{StaSep07}.} from those in \cite{StaSep07}, hold:

H1. There exists a constant $\gamma$ such that $\phi'(y)\le \gamma, \ \forall y\in R.$

H2.  The triplet $(A+\gamma^* BC, B,C)$ is {\it strictly positive real} \cite{khalil}; that is, there exits $P=P^T>0$ such that:
\begin{eqnarray}\label{kyp}
P(A+\gamma^*BC)+(A+\gamma^*BC)^TP < 0 \\
PB=C^T.\label{kypeq}
\end{eqnarray}

H3.  $L1_N=0$ and the largest number, $\lambda_2$, such that (\ref{luca}) holds, satisfies:
\begin{equation}\label{gammaz}
\lambda_2>\gamma-\gamma^*.
\end{equation}

To apply Theorem \ref{odethm}, note that system (\ref{stan1})-(\ref{stan3}) is  of the form (\ref{compa})-(\ref{Fnot}), with
\begin{equation}
f(x)=Ax+B\phi(Cx) \quad D=BC.
\end{equation}
From H1, we conclude that the Jacobian $J(x)$ is as in (\ref{membership}), with
\begin{equation}\label{z1s1}
Z_1=A+\gamma BC \quad S_1=-BC.
\end{equation}
Noting from (\ref{kypeq}) that
\begin{equation}\label{important}
PD=PBC=C^TC\ge 0
\end{equation}
and using (\ref{gammaz}), we obtain:
\begin{equation}\label{trick}
P(Z_1-\lambda_2D)+(Z_1-\lambda_2D)^TP\le P(Z_1-(\gamma-\gamma^*)D)+(Z_1-(\gamma-\gamma^*)D)^TP.
\end{equation}
Substituting $Z_1$ from (\ref{z1s1}) in the right-hand side of (\ref{trick}) and using (\ref{kyp}), we conclude that condition (\ref{hull1}) of Theorem \ref{LMIlem} holds.  Likewise, (\ref{hull2}) holds because $PS_1=-PBC=-C^TC\le 0$, and Theorem \ref{LMIlem} verifies condition (\ref{diagstab}) of Theorem \ref{odethm} on $\mathcal{X}=\mathbb{R}^n$. Finally, noting from (\ref{important}) that (\ref{diag}) also holds, and that $PD$ is symmetric as stipulated in Remark \ref{directed}, we conclude (\ref{sync}) for all forward-complete\footnote{\cite{StaSep07} indeed argues boundedness for the solutions of (\ref{stan1})-(\ref{stan3}), using ideas from \cite{ArcTee02}.} trajectories. \hfill $\Box$
\end{example}

\section*{Appendix:  Proofs for Theorem \ref{LMIlem} and Theorem \ref{dyad}}\label{app}

\noindent
{\bf Proof of Theorem \ref{LMIlem}:} From (\ref{membership}), for every $x\in \mathcal{X}$, there exist parameters $\theta_1,\cdots,\theta_q,\omega_1,\cdots,\omega_m\ge 0$, $\theta_1+\cdots+\theta_q=1$, such that:
\begin{equation}
J(x)=\theta_1Z_1+\cdots+\theta_qZ_q+\omega_1S_1+\cdots+\omega_mS_m.
\end{equation}
Because $\theta_1+\cdots+\theta_q=1$, we write:
\begin{equation}\label{above}
J(x)-\lambda_2D =\theta_1(Z_1-\lambda_2D)+\cdots+\theta_q(Z_q-\lambda_2D)+\omega_1S_1+\cdots+\omega_mS_m,
\end{equation}
from which it follows that a matrix $P$ satisfying (\ref{hull1})-(\ref{hull2}) also satisfies (\ref{diagstab}).  To prove the converse, we note from the surjectivity assumption that, for any set of parameters $\theta_1,\cdots,\theta_q,\omega_1,\cdots,\omega_m\ge 0$, $\theta_1+\cdots+\theta_q=1$, there exits $x\in \mathcal{X}$ for which (\ref{above}) holds.  To see that (\ref{diagstab}) implies (\ref{hull1}), pick $\theta_k=1$, $\theta_i=0 \ i\neq k$, and $\omega_i=0 \ i=1,\cdots,m$ in (\ref{above}). To see that (\ref{diagstab}) implies (\ref{hull2}), assume, to the contrary, that (\ref{diagstab}) holds, but (\ref{hull2}) fails for some $k$, which means that there exists $\zeta \in \mathbb{R}^n$ such that
\begin{equation}
\zeta^T(PS_k+S_k^TP)\zeta>0.
\end{equation}
Then, pick $\omega_i=0\ i\neq k$ and note from (\ref{above}) that the left-hand side of (\ref{diagstab}) is equal to:
\begin{equation}
\omega_k(PS_k+S_k^TP)+\sum_{k=1}^q \theta_k[P(Z_k-\lambda_2D)+(Z_k-\lambda_2D)^TP].
\end{equation}
Because $0\le \theta_k\le 1$, choosing $\omega_k>0$ large enough ensures that
\begin{equation}
\zeta^T\left\{\omega_k(PS_k+S_k^TP)+\sum_{k=1}^q \theta_k[P(Z_k-\lambda_2D)+(Z_k-\lambda_2D)^TP]\right\}\zeta>0,
\end{equation}
which contradicts (\ref{diagstab}). \hfill $\Box$
\bigskip

\noindent
{\bf Proof of Theorem \ref{dyad}:}  We rewrite (\ref{blockdia})-(\ref{composite}) as
\begin{equation}\label{dortbucuk}
\left[\begin{array}{cc} P(A_0-\lambda_2D)+(A_0-\lambda_2D)^TP & PB+CQ \\ QC^T+B^TP & -2Q\end{array}\right]<0,
\end{equation}
where $Q=diag\{q_1,\cdots,q_\ell\}$, and make use of the following lemma, proven separately below:

\begin{lemma}\label{S}
If there exists a diagonal $\ell\times \ell$ matrix $Q>0$ satisfying (\ref{dortbucuk}),
then
\begin{equation}\label{sifir}
P(A_0+\gamma_1 B_1C_1^T+\cdots+\gamma_\ell B_\ell C_\ell^T-\lambda_2D)+(A_0+\gamma_1 B_1C_1^T+\cdots+\gamma_\ell B_\ell C_\ell^T-\lambda_2D)^TP<0 \quad \forall \gamma_i \in [0,1].
\end{equation}
When $\ell=1$, the converse is also true; that is, if (\ref{sifir}) holds for every $\gamma_1\in [0,1]$, then (\ref{dortbucuk}) holds for some constant $Q>0$.
\end{lemma}

\noindent
 To conclude the first statement of Theorem \ref{dyad}, we note that (\ref{sifir}) implies (\ref{diagstab}) with $\epsilon>0$.  To prove the second statement, we note from the surjectivity assumption that (\ref{diagstab}) with $\epsilon>0$ implies (\ref{sifir}).  Since $\ell=1$, we apply the converse statement in Lemma \ref{S} and conclude that (\ref{dortbucuk}) holds for some constant $Q>0$; that is, (\ref{composite}) holds with $\ell=1$ and $\mathcal{P}$ as in (\ref{blockdia}).
\hfill $\Box$

\bigskip

\noindent
{\bf Proof of Lemma \ref{S}:} Defining $\Gamma$ to be a diagonal matrix with entries $\gamma_i \in [0,1]$, we rewrite (\ref{sifir}) as
\begin{equation}\label{bir}
x^T[P(A_0-\lambda_2D+B\Gamma C^T)+(A_0-\lambda_2D+B\Gamma C^T)^TP]x<0 \quad \forall x\neq 0.
\end{equation}
We then define the new variable
\begin{equation}\label{defy}
y:=\Gamma C^Tx,
\end{equation}
and rewrite (\ref{bir}) as:
\begin{equation}\label{iki}
[x^T \ y^T]\left[\begin{array}{cc} P(A_0-\lambda_2D)+(A_0-\lambda_2D)^TP & PB \\ B^TP & 0\end{array}\right]\left[\begin{array}{c} x \\ y \end{array}\right]<0.
\end{equation}
Next, we note from (\ref{defy}) with $\gamma_i \in [0,1]$ that $y_i$ and $C_i^Tx$ are constrained by:
\begin{equation}
y_i(C_i^Tx)=\frac{1}{\gamma_i}y_i^2 \ge y_i^2,
\end{equation}
which means that:
\begin{equation}\label{uc}
[x^T \ y^T]\left[\begin{array}{cc} 0 & C_ie_i^T \\ e_iC_i^T & -2e_ie_i^T\end{array}\right]\left[\begin{array}{c} x \\ y \end{array}\right]\ge 0  \quad i=1,\cdots,\ell,
\end{equation}
where $e_i$ is the $i$th unit vector in $\mathbb{R}^\ell$. Thus, (\ref{bir}) is equivalent to the statement that (\ref{iki}) holds for all $x\neq 0$, $y\neq 0$, satisfying (\ref{uc}) $i=1,\cdots,\ell$.

We now invoke the {\it S-procedure} \cite{boyd} which states that, for symmetric matrices $T_0,T_1,\cdots,T_\ell$,
\begin{equation}\label{u}
\zeta^TT_0\zeta <0 \quad \mbox{for all $\zeta \neq 0$ satisfying} \quad \zeta^TT_i\zeta\ge 0 \quad i=1,\cdots,\ell
\end{equation}
if there exist $q_1>0,\cdots,q_\ell>0$ such that
\begin{equation}\label{d}
T_0+q_1 T_1+\cdots+q_\ell T_\ell<0.
\end{equation}
Because the matrices in (\ref{iki}) and (\ref{uc}) play the roles of $T_0$ and $T_i$, $i=1,\cdots,\ell,$ in the S-procedure, we conclude that (\ref{iki}) holds for all $x\neq 0$, $y\neq 0$, satisfying (\ref{uc}) if
\begin{equation}\label{dort}
\left[\begin{array}{cc} P(A_0-\lambda_2D)+(A_0-\lambda_2D)^TP & PB \\ B^TP & 0\end{array}\right]+\left[\begin{array}{cc} 0 & CQ \\ QC^T & -2Q\end{array}\right]<0
\end{equation}
for some diagonal $Q\ge 0$.  Finally, we note that (\ref{dort}) requires $Q>0$ because, if $Q$ contains zero diagonal entries,  then the matrix in (\ref{dort}) also contains zero diagonal entries and, thus, cannot be negative definite.  This concludes the proof of the first statement of the lemma, because inequality (\ref{dort}) is identical to (\ref{dortbucuk}).

To prove the converse statement, we recall that, when $\ell=1$, the S-procedure also states that (\ref{u})  implies (\ref{d}) for some $q_1\ge 0$, provided there exists $\zeta_0$ such that $\zeta_0^TT_1\zeta_0>0$. When $\ell=1$, $T_1$ defined in (\ref{uc}) has the form:
\begin{equation}
T_1=\left[\begin{array}{cc} 0 & C\\ C^T & -2\end{array}\right],
\end{equation}
which means that $\zeta_0^TT_1\zeta_0>0$ indeed holds with the choice $\zeta_0=[x_0^T \ \frac{1}{2} (C^Tx_0)]^T$, where $x_0$ is such that $C^Tx_0\neq 0$.  Because $T_0$ is as defined in (\ref{iki}), we conclude from the S-procedure that if (\ref{sifir}) holds for every $\gamma_1\in [0,1]$, then (\ref{dort}) holds for some constant $Q\ge 0$.  Recalling that (\ref{dort}) cannot hold if $Q=0$ and that (\ref{dort}) is identical to (\ref{dortbucuk}), we conclude that (\ref{dortbucuk}) must be true for some $Q>0$.
\hfill $\Box$


\end{document}